\documentclass[12pt,twoside]{amsbook}
\usepackage{graphicx,amscd,amssymb,verbatim}

\theoremstyle{plain}
\newtheorem{Thm}{Theorem}[section]
\newtheorem{Cor}[Thm]{Corollary}

\newtheorem{Lem}[Thm]{Lemma}

\theoremstyle{definition}
\newtheorem{Def}[Thm]{Definition}

\theoremstyle{remark}
\newtheorem{Remark}[Thm]{Remark}

\newtheorem{Example}[Thm]{Example}

\numberwithin{equation}{section}
\numberwithin{section}{chapter}

\numberwithin{figure}{chapter}

\DeclareMathOperator{\rank}{rank}

\def\vs#1{\vskip .#1 cm}

\def\picture #1 by #2 (#3){
  \vbox to #2 {
     \hrule width #1 height 0pt depth 0pt
     \vfill
     \special{picture #3} 
  }
}

\def\scaledpicture #1 by #2 (#3 scaled #4){{
   \dimen0=#1 \dimen1=#2
   \divide\dimen0 by 1000 \multiply\dimen0 by #4
   \divide\dimen1 by 1000 \multiply\dimen1 by #4  
   \picture\dimen0 by dimen1 (#3 scaled #4)}
}

\rm
\begin{document}


\newpage
\pagenumbering{roman}
\begin{center}

\vskip1cm

{\Large Generalizations of two-stack-sortable permutations}

\vskip2cm
A Dissertation
\vskip0.6cm
{Presented to}
\vskip0.6cm
{The Faculty of the Graduate School of Arts and Sciences}
\vskip0.6cm
{Brandeis University}
\vskip0.6cm
{Department of Mathematics}
\vskip0.6cm
{Professor Ira Gessel, Advisor}
\vskip0.6cm
{In Partial Fulfillment}
\vskip0.2cm
{of the Requirements for the Degree}
\vskip0.2cm
{Doctor of Philosophy}
\vskip0.6cm
{by}
\vskip0.6cm
{Dapeng Xu}
\vskip0.6cm
{September, 2002}
\end{center}

\newpage
\pagenumbering{roman}
\setcounter{page}{3}
\begin{center}
{\huge Acknowledgment}
\end{center}
\vskip2cm

\indent My deepest gratitude goes to my advisor, Ira Gessel, without whom this work would not be possible; I
would like to thank him for his generosity in sharing his insight and time with me, and his constant patience
and encouragement. It was and will always be my pleasure to work with him.

I would like to thank Daniel Ruberman, Kiyoshi Igusa, Alan Mayer and Mark Adler for teaching me a lot of
mathematics. I would like to thank Susan Parker for helping me to be a better teacher.

I'd also like to thank the Mathematics Department faculty and  staff for making my life here more enjoyable.

\vskip1cm

\newpage
\vskip 3cm
\begin{center}
{\Large ABSTRACT}
\vskip1cm
{\bf Generalizations of two-stack-sortable permutations}
\vskip1cm
A dissertation presented to the faculty of 

the Graduate School of Arts and Sciences of 

Brandeis University, Waltham, Massachusetts
\vskip1cm
by Dapeng Xu
\end{center}
\vskip1cm
In this thesis, we apply the stack sorting operator to $r$-permutations and 
construct the functional equation for the generating
function of two-stack-sortable $k$-tuple $r$-permutations counted by descents by using a factorization
similar to Zeilberger's. We solve the functional equation and give explicit formulas for
the number of two-stack-sortable $r$-permutations.



\newpage
\pagenumbering{arabic}
\setcounter{page}{1}
\tableofcontents

\chapter{Introduction}

The operation of stack sorting was first studied by Knuth [\ref{Knuth}, p.~ 239].
He described the operation by a railroad cars swapping algorithm. Later West described
this operation in terms of a simple card game [\ref{West}]. 
 
 The operation uses a stack and can be described as follows. Let $\pi$ be
a word with distinct letters in the alphabet $\{1,2,3,\dots\}$. We call $\pi$ 
the input of the operation. At the $i$th step, we compare the $i$th
letter $a$ of the word $\pi$ with the top letter on the stack (if
any). If $a$ is smaller or the stack is empty, we put $a$ on the top of the stack.
If
$a$ is bigger, we move the top letter on the stack to the output and then again
compare $a$ with the new top letter on the stack (if any). We repeat this
until we put $a$ onto the stack. When the input is empty, we move all the letters
 on the stack (if any) to the output in order from top to bottom. 

West also gave a recursive
definition of the stack sorting operation.
[\ref{West}, \ref{West2}].
\begin{Def}
Let $\pi=a_1a_2\cdots a_n$ be a word in the alphabet ${\mathbb
P}=\{1,2,\dots\}$, having all its letters distinct. If $n=0$, let
$S(\pi)$ be the empty word. Otherwise, let $S(\pi)$ be obtained by
permuting the letters of $\pi$ as follows: if $m=\max (a_1, a_2, \dots,
a_n)$ and $\pi=\pi_l m \pi_r$, then 
$$ S(\pi)=S(\pi_l)S(\pi_r)m.$$
\end{Def}

It is easy to see that this recursive definition of the stack sorting
operation is equivalent to the description in terms of stacks. For a given word
$\pi=a_1a_2\cdots a_n$ in the alphabet ${\mathbb
P}=\{1,2,\dots\}$, having all its letters distinct, we can always
write it as $\pi=\pi_l m \pi_r$ where $m=\max (a_1, a_2, \dots,
a_n)$. According to the description, since $m$ is bigger than any letter
in $\pi_l$, all the letters in $\pi_l$ will be moved to the output before  we can put
$m$ onto the stack. Also since $m$ is bigger than any letter in $\pi_r$,
$m$ will not be removed from the stack until all the letters in $\pi_r$
go through the stack. That is exactly what the recursive definition says.

The problem of stack sorting was then generalized and reseached in a
number of ways. Among these, the problems of enumerating $t$-stack-sortable permutations
interest most people. We say that a permutation $\pi$ of $[n]$ is a $t$-stack-sortable
permutation if 
$S^t(\pi)$ is the identity permutation, and that $\pi$ is exactly $t$-stack-sortable if
$\pi$ is
$t$-stack-sortable but not $(t-1)$-stack-sortable. Knuth [\ref{Knuth}, p.~ 239] proved
that the number of 1-stack-sortable permutations of $[n]$ is the Catalan number
$C_n=\binom{2n}{n} /(n+1)$. West [\ref{West2}, \ref{West}] proved that all
permutations of $[n]$ are $(n-1)$-stack-sortable, that $(n-2)!$ are exactly
$(n-1)$-stack-sortable, and that    $\frac{7}{2}(n-2)!+(n-3)!$ are exactly
$(n-2)$-stack-sortable. 

Characterizations of 1-stack-sortable permutations and 
2-stack-sortable permutations were also given using the following notion of pattern avoidance. 
 Let $ \pi_1=a_1a_2\cdots a_n $ be a permutation of $[n]$
and let $ \pi_2=b_1b_2 \cdots b_k$ be a permutation of $[k]$. We say that $\pi_1$
contains a type $\pi_2$ subsequence if there exists $1 \le i_{b_1} < i_{b_2} < \cdots
< i_{b_k} \le n$ such that $a_{i_1} < a_{i_2} < \cdots < a_{i_k}$. We say $\pi_1$
avoids $\pi_2$ if $\pi_1$ contains no $\pi_2$-subsequence. 
For example, a
permutation $\pi$ avoids 231 if it cannot be written as $\cdots a \cdots b \cdots c
\cdots$ so that $ c <a<b$. 

Tarjan [\ref{tarjan}] proved that  a permutation is 1-stack-sortable if and only if it avoids the pattern
231. West [\ref{West2}, \ref{West}] proved that a permutation fails to be
2-stack-sortable if it contains a subsequence of type 2341 or a
subsequence of type 3241 which is not part of a subsequence of type 35241
and it is two-stack-sortable if it contains no such subsequence. West also
conjectured the number of 2-stack-sortable permutations of $[n]$ to be [\ref{West},
\ref{West2}]
$$\frac {2(3n)!}{(n+1)!\,(2n+1)!}.$$
The conjecture was first proved by D. Zeilberger
[\ref{Zeil}]. Later, two bijections between two-stack-sortable permutations
and non-separable planar graphs were given by S. Dulucq, S. Gire and  O.
Guibert [\ref{DGG}]; I. P. Goulden and J. West [\ref{GW}]. More contributions to
this problem were given by Mikl{\'o}s B{\'o}na, Mireille
Bousquet-M{\'e}lou, Leopold Travis and others [\ref{Bona1}, \ref{Bona2},
\ref{melou1}, \ref{melou2}, \ref{melou3}, \ref{Travis}].

\section{The generalization}

In this thesis, we generalize the ordinary permutations to $r$-permutations [\ref{PC},
\ref{PC2}, \ref{park1}, \ref{park2}, \ref{park3}] and then enumerate
the number of 2-stack-sortable $r$-permutations under the generalized stack-sorting
operation on $r$-permutations.

\begin{Def}
If a permutation $a_1 a_2 \cdots a_{rn}$ of $\{1^r, 2^r,
\dots , n^r \} $ satisfies the condition that if $i<j<k$ and $a_i =a_k$
then
$a_j
\le a_i$, we call it an {\it $r$-permutation} of $[n]$.
\end{Def}

\begin{Def} \label{stacksortingoperation}
Let $\pi$ be an
$r$-permutation of $[n]$. For $n>0$,  we can write $\pi$ as $\pi=\alpha_1 n
\alpha_2 n
\cdots n \alpha_{r+1}$. The stack-sorting operation $S$ on
$r$-permutations is defined by
$$ S(\pi)=S(\alpha_1)S(\alpha_2)\cdots S(\alpha_{r+1})n.$$
When $n=0$, i.e., when $\pi$ is an empty permutation, we define $S(\pi)$ to be $\emptyset$. 
Note that $S(\pi)$ is
an ordinary permutation.
\end{Def}

\begin{Def}
Given an $r$-permutation $\pi$ and a letter $a$ in $\pi$, we call $a$ a type $i$
descent
$(i=1,\dots,r)$ if the $i$th occurrence of $a$ is immediately followed by a smaller
letter, and we call $a$ a type 0 descent if the first occurrence of $a$ immediately
follows a smaller letter. We denote the set of type $i$ descents of an $r$-permutation
$\pi$ by
$\pi^{(i)}$ 
$(i=0,\dots,r)$.
\end{Def}

\begin{Remark}
Notice that the $r$-permutations are usually defined the other way, i.e., a permutation $a_1a_2\cdots a_{rn}$
of $\{1^r, 2^r,
\dots , n^r \} $ is an $r$-permutation if it satisfies the condition that if  $i<j<k$
and $a_i=a_k$, then
$a_j
\ge a_i$. The reason for the change is that we expect $S(\pi)$ to be an ordinary
permutation for an
$r$-permutation $\pi$, while the traditional definition of $r$-permutation can not give
us that. For example, 123321 is a 2-permutation under the traditional definition, but
$S(123321)=121233$. This behavior is very different from what we study here.
\end{Remark}

\begin{Remark}
 Note that a type 0
descent is actually an ascent. 
\end{Remark}
\begin{Remark}
Often, the last letter (or, the last position) of a permutation is considered to be a
descent but not here.
\end{Remark}

\begin{Example} \label{ex1}
When $r=3$, $n=6$,
\[
\pi=544453222335611166
\]
is a $3$-permutation of $\{1,2,\dots,6\}$. Then
\begin{align*}
S(\pi)=&S(544453222335)S(111)6 \\
=&S(444)S(322233)516 \\
=&4S(222)3516\\
=&423516.
\end{align*}
The sets  of type $i$ descents
$(i=0,\dots,3)$ of $\pi$ are 
$\pi^{(0)}=\{6\}$, $\pi^{(1)}=\{3,5,6\}$, $\pi^{(2)}=\{5\}$ and $\pi^{(3)}=\emptyset$.
\end{Example}

\begin{figure} \label{treerep}
\centering
\includegraphics[scale=.7]{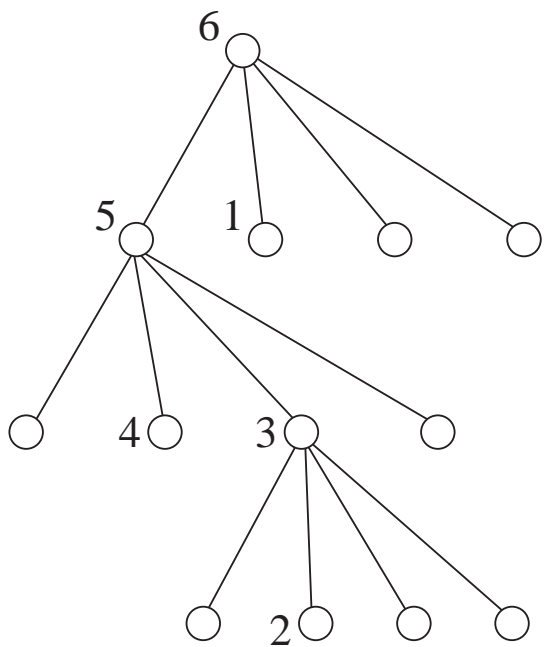}
\caption{ Tree representation of $\pi$}
\end{figure}

The $r$-permutations can be represented as $(r+1)$-ary  decreasing trees. If we
have an $r$-permutation $\pi$ of $[n]$ and $\pi=\alpha_1 n \alpha_2 n \cdots n
\alpha_{r+1}$, to get the tree representation of $\pi$, we set $n$ to be the root
of the tree and recursively set the tree representation of $\alpha_i$ ($i=1,2,\dots,
r+1$) to be the $i$th child of the root.
  If $\alpha_i$ is nonempty, then
$\pi$ has an $i$th child and $n$ is a type $i-1$ descent. Therefore it is clear that each $i$th 
child corresponds to a type $i-1$
descent. For example, the
3-permutation above can be represented as in Figure
\ref{treerep}.

There are three kinds of traversals of a tree: preorder, inorder, and
postorder [\ref{tamassia}, p.~ 243; \ref{Knuth}, p.~ 315]. In a preorder traversal of a
tree, the root is visited first and then the subtrees rooted at its children are
traversed recursively. A postorder traversal recursively traverses the subtrees
rooted at the children of the root first, and then visits the root. 
An inorder traversal, in our case, recursively traverses the subtree
rooted at the first child of the root, then visits the root, then
recursively traverses the subtree rooted at the second child of the root,
then visits the root, and so on. Hence,
 if we read the tree in Figure \ref{treerep} in inorder, we get back the
original
$r$-permutation
$\pi$. If we read it in postorder, we get an ordinary permutation and it  
corresponds to $S(\pi)$.  

Because any $(r+1)$-ary tree on $n$ vertices has $n-1$ children, we have
\begin{Lem} \label{tree}
Every $r$-permutation of $[n]$ has a total of  $n-1$ descents of types $0,1,\dots,r$.
\end{Lem}

\section{Stack-sortable $r$-permutations}

Now let us first consider the enumeration of 1-stack-sortable $r$-permutations (which we
call here stack-sortable $r$-permutations) with all types of descents.

 We weight  any type
$i$ descent by $x_i$. Therefore the weight of the $r$-permutation $\pi$ in  Example 
\ref{ex1} is
$x_0x_1^3x_2$.
Because of Lemma \ref{tree},  we do not need a parameter to keep track of the number of
different letters in an $r$-permutation.

 For any unlabeled $(r+1)$-ary tree, there is one and only one way
to label the nodes to make it an $(r+1)$-ary decreasing tree such that if
we read it in postorder, we get an identity ordinary permutation. That gives a
bijection between stack-sortable $r$-permutations and unlabeled $(r+1)$-ary trees.
Therefore, 
if we weight an unlabeled $(r+1)$-ary tree by weighting any $i$th child by
$x_{i-1}$, then the total weight of stack-sortable $r$-permutations is the total
weight of unlabeled
$(r+1)$-ary trees. Let 
$A(x)$ be the weight of all  stack-sortable
$r$-permutations. Then we have
$$
A(x)=\prod_{i=0}^r (1+x_iA(x))
$$
where $x=(x_0,x_1,\dots, x_r)$. 

To solve this functional equation, we
use  Lagrange inversion [\ref{Goulden}, p.~ 21; \ref{Ira}]. 

\begin{Lem} \label{lagrange}
{\mbox (Multivariable Lagrange Inversion)} Let
$f(\lambda)
\in {\mathbb R}[[\lambda]]$ and
$\phi_1(\lambda),\dots,$ $\phi_m(\lambda) \in {\mathbb R}[[\lambda]]$ where
$\lambda=(\lambda_1,\dots,\lambda_m)$. Suppose that $w_i=t_i\phi_i(w)$
for $i=1,\dots,m$, where $w=(w_1,\dots,w_m)$. Let
$\phi=(\phi_1,\dots,\phi_m)$, $t=(t_1,\dots,t_m)$ and $k=(k_1,\dots,k_m)$. Then
$$f(w(t))=\sum_k t^k [\lambda^k] \left\{ f(\lambda)\phi^k(\lambda)
\left| {\delta_{ij}- \frac {\lambda_j}{\phi_i(\lambda)} \frac {\partial
\phi_i(\lambda)}{\partial \lambda_j}} \right|\right\},$$
where $[\lambda^k]f(\lambda)$ is the coefficient of
$\lambda_1^{k_1}\cdots\lambda_m^{k_m}$ of the formal power series
$f(\lambda)$.
\end{Lem}

In our case, we let $w_i=x_iA(x)$ and
$\phi_i(\lambda)=\prod_{j=0}^r (1+\lambda_j)$ for
$i=0,\dots,r$. Then $w_i$ satisfy the condition
that $w_i=x_i\phi_i(w)$. Also let
$f(\lambda)=\prod_{j=0}^r(1+\lambda_j)$. Then $f(w)=A(x)$.
Using multivariable Lagrange inversion, we get 
\begin{align*}
A(x)=&f(w) \\
=& \sum_k x^k [\lambda^k] \left\{ \left(\prod_{i=0}^r (1+\lambda_i)
\right)^{k_0+\cdots+k_r+1}
\left| {\delta_{ij}- \frac {\lambda_j}{1+\lambda_j} } \right|\right\}.
\end{align*}
Now we introduce some notation. Let
\begin{equation}
{ e_k(x)=\sum_{0 \le i_1 < i_2 < \cdots <i_k
\le r} x_{i_1} x_{i_2} \cdots x_{i_k} }
\end{equation}     
be the $k$th elementary symmetric function of $x=(x_0, x_1,\dots,x_r)$. 
Define 
\begin{align} \label{defE1}
E(x,u)=&E(x_0,x_1,\dots,x_r;u)\\
=&\prod_{i=0}^{r} (1+x_iu) \notag\\
=&\sum_{i=0}^{r}e_i(x)u^i \notag
\end{align}
and let 
\begin{equation} \label{E(x)}
E(x)=E(x,1).
\end{equation} 
Also we introduce the following lemma to evaluate the determinant.

\begin{Lem} \label{determinant}
Let $M(x)=\left|\delta_{ij}+(\delta_{ij}-1)x_j\right|$,  where $0 \le i , j \le r$, be
the determinant of the 
$(r+1) \times (r+1)$ matrix. Then
$$ M(x)=E(x)\left(1-\sum_{i=0}^{r} \frac{x_i}{1+x_i}\right). $$
\end{Lem}
\begin{proof}
\[
M(x)= \left| \begin{array}{cccccc}
         1  &-x_1& -x_2& \cdots & -x_{r-1} &-x_r \\
         -x_0 &1  &-x_2 &\cdots & -x_{r-1} &-x_r \\
	        \cdots \cdots\\
         -x_0 & -x_1 &-x_2 &\cdots &1  &-x_r \\
         -x_0 &-x_1 &-x_2 &\cdots& -x_{r-1} &1
             \end{array} 
      \right| .
\]
Subtracting the last row from every other row, we get 
\[
 M(x) = \left| \begin{array}{cccccc}
         1+x_0 & 0 & 0 &\cdots &0 & -(1+x_r) \\
         0   & 1+x_1 & 0 &\cdots &0 & -(1+x_r) \\
         0 & 0 & 1+x_2 & \cdots &0 & -(1+x_r) \\
         \cdots \cdots \\
         0 & 0 & 0 & \cdots &1+x_{r-1} & -(1+x_r) \\
         -x_0 & -x_1 & -x_2 & \cdots &-x_{r-1} & 1
              \end{array}
           \right| .
\]    
If we let $x^{\prime}=(x_1,x_2,\dots,x_r)$, then by expanding the first column, we get
\begin{align*}
 M(x) =&(1+x_0)
         \left| \begin{array}{ccccc}
         1+x_1 & 0  &\cdots  & -(1+x_r) \\
         0   & 1+x_2  &\cdots  & -(1+x_r) \\
         \cdots \cdots \\
         -x_1 & -x_2  & \cdots  & 1
              \end{array}
           \right|   \\
								&-(-1)^r x_0
								\left| \begin{array}{ccccc}
          0 & 0 &\cdots & 0 &-(1+x_r) \\
          1+x_1 & 0 &\cdots & 0 &-(1+x_r) \\
         \cdots \cdots \\
         0 & 0 & \cdots  & 1+x_{r-1}& 1
              \end{array}
           \right|    \\
=&(1+x_0)M(x^{\prime}) -(-1)^r x_0 (-1)^r \prod_{i=1}^r(1+x_i) \\
 =&(1+x_0)  M(x^{\prime}) -x_0 \prod_{i=1}^r(1+x_i).
\end{align*}

Now we use induction. It is easy to see that the lemma is true when
$r=0$. Suppose the lemma is true for
$M(x^{\prime})$, i.e., 
\[
M(x^{\prime})=E(x^{\prime})\left(1- \sum_{i=1}^{r} \frac{x_i}{1+x_i} \right).
\]
Then 
\begin{align*}
M(x)=&(1+x_0)E(x^{\prime})\left(1-\sum_{i=1}^{r} \frac{x_i}{1+x_i}
\right)-x_0 \prod_{i=1}^r(1+x_i) \\
=&E(x)\left(1-\sum_{i=1}^{r} \frac{x_i}{1+x_i}
\right)-E(x) \frac{x_0}{1+x_0}\\
=&E(x)\left(1-\sum_{i=0}^{r} \frac{x_i}{1+x_i} \right).
\end{align*}
The lemma is proved.
\end{proof}

With Lemma \ref{determinant}, it follows that 
\begin{align*}
A(x)=& \sum_k x^k [\lambda^k] \left\{ \left(\prod_{i=0}^r (1+\lambda_i)
\right)^{k_0+\cdots+k_r+1}
\left(\sum_{i=0}^r \frac {1}{1+\lambda_i} -r \right) \right\} \\
=& \sum_k x^k \left\{ \sum_i \binom{n}{k_1} \cdots \binom{n-1}{k_i}
\cdots \binom {n}{k_r} - r \prod_i \binom{n}{k_i} \right\} \\
=& \sum_k x^k \left\{ \left( \sum_i \frac{n-k_i}{n} -r \right) \prod_i 
\binom{n}{k_i} \right\} \\
=& \sum_k x^k  \frac {1}{n} \prod_{i=0}^r \binom {n}{k_i}
\end{align*}
where $n=1+k_0+\cdots+k_r$. Therefore, 
\begin{Thm}
The number of stack-sortable
$r$-permutations with $k_i$ descents of type $i$ is  
$$\frac {1}{n} \prod_{i=0}^r \binom {n}{k_i}$$ 
where $n=1+k_0+\cdots+k_r$.
\end{Thm}

\begin{Remark}
This functional equation can also be solved by
one variable Lagrange inversion by introducing a new variable,
i.e., letting 
$$A(z)=z\prod_{i=0}^r (1+x_iA(z)).$$
\end{Remark}

\begin{Remark}
These numbers also come up in counting noncrossing partitions
[\ref{stan}, \ref{simion}, \ref{edelman}]. 
\end{Remark}

\begin{Remark}
When $r=1$, these numbers are Narayana numbers [\ref{sulanke},
\ref{hwang}].
\end{Remark}
\vs 5

\chapter{Two-stack-sortable $r$-permutations}
\section{The functional equations} \label{mainfuneq}

Now let us consider the case of two-stack-sortable $r$-permutations. An $r$-permutation
$\pi$ of $[n]$ is two-stack-sortable if $S^2(\pi)$ is an identity permutation.
Let
$\pi=\alpha_1 n
\alpha_2 n
\cdots n \alpha_{r+1}$ where $n$ is the largest element in $\pi$. Then
\begin{align*}
S(\pi)=&S(\alpha_1)S(\alpha_2)\cdots S(\alpha_{r+1})n \\
S^2(\pi)=&S(S(\alpha_1)S(\alpha_2)\cdots S(\alpha_{r+1}))n
\end{align*}

Denote the identity permutation by $I$. (Here we abuse the notation by letting $I$  be
the identity permutation of $[n]$ for any positive integer $n$.) One can notice that if
$S^2(\pi)=I$ and if
$m$ is the largest number in 
$\alpha_1  \alpha_2 \cdots  \alpha_{r+1}$, then $m$ can occur in only 
one of the
$\alpha_i$ because if $m$ appeared in two of the $\alpha_i$, the $r$-permutation
condition for $\pi$ would be violated. Thus
$(\alpha_1,
\alpha_2,
\dots,
\alpha_{r+1})$ is an $(r+1)$-tuples of $r$-permutations satisfying the
condition that 
$$S(S(\alpha_1)S(\alpha_2)\cdots S(\alpha_{r+1}))=I.$$
This suggests that we study the subject defined as follows.
\begin{Def}
 $(\alpha_1, \alpha_2, \dots,
\alpha_k)$ is a  $k$-tuple $r$-permutation of $[n]$ if $\alpha_1
\alpha_2 \cdots \alpha_k$ is an $r$-permutation of $[n]$ and any letter
appears in only one of  the $\alpha_i$. Also we call $\alpha_i$ the $i$th component of the
$k$-tuple $r$-permutation.
\end{Def}
We call a $k$-tuple $r$-permutation $(\alpha_1, \alpha_2, \dots,
\alpha_k)$ a two-stack-sortable $k$-tuple $r$-permutation if it satisfies the condition
that 
$$S(S(\alpha_1)S(\alpha_2)\cdots S(\alpha_{k}))=I.$$

From the definition, we know that each component
$\alpha_i$ of a $k$-tuple $r$-permutation is itself an
$r$-permutation so that we can interpret a $k$-tuple $r$-permutation of $[n]$ as a
forest of $k$ $(r+1)$-ary decreasing trees.  The set of the descents of type $i$ of
$k$-tuple
$r$-permutations
$\pi=(\alpha_1,
\dots,
\alpha_k)$ is the union of the sets of descents of type $i$ of $\alpha_j$ $(j=1,\dots
k)$, i.e.
\[ \displaystyle
\pi^{(i)}=\bigcup_{j=1}^k \alpha_j^{(i)}.
\]

\begin{figure} \label{treerep2}
\centering
\includegraphics[scale=.7]{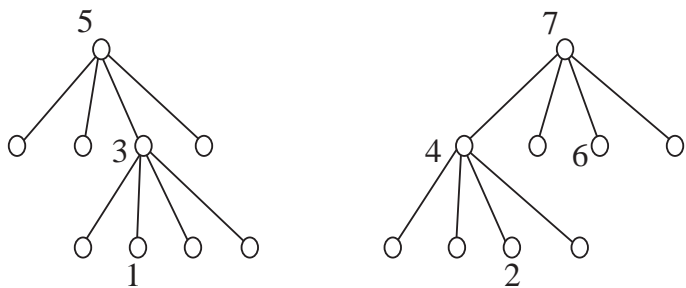}
\caption{ Forest representation of $\pi$}
\end{figure}

For example,  let $\pi=(\alpha_1, \alpha_2)$ where  $\alpha_1=553111335$ and
$\alpha_2=442224776667$. Then
$(\alpha_1, \alpha_2)$ is a 2-tuple 3-permutation of $\{1,2,\dots
,7\}$ and it can be represented as in Figure \ref{treerep2}.
Since 
\begin{align*}
S(S(\alpha_1)S(\alpha_2))=& S(1352467) \\
=&1324567,
\end{align*}
$\pi$ is not two-stack-sortable. The sets of descents of $\pi$ are $\pi^{(0)}=\{7\}$, 
$\pi^{(1)}=\{3\}$, $\pi^{(2)}=\{4,5,7\}$, $\pi^{(3)}=\emptyset$.  

Let $g_n^{(k)}$ be the sum of the weights of two-stack-sortable $k$-tuple
$r$-permutations of $n$ different letters  
{\em such that every component is nonempty}. We need every component to be nonempty
so that we can keep track of every type of descent. Suppose $(\alpha_1, \alpha_2,
\dots,
\alpha_{k})$ is a two-stack-sortable $k$-tuple $r$-permutation of $[n]$ in which every
component is nonempty and the largest element
 $n$ appears in $\alpha_l$. Let $\alpha_l=\beta_1 n \beta_2 n
\cdots n
\beta_{r+1}$. Then

\begin{align} \label{decomposition}
\\
I=&S(S(\alpha_1)S(\alpha_2)\cdots S(\alpha_{k})) \notag \\ 
 =&S(S(\alpha_1)S(\alpha_2)\cdots S(\alpha_{l-1}) S(\beta_1)  S(\beta_2) 
\cdots S(\beta_{r+1}) \, n \, S(\alpha_{l+1}) \cdots S(\alpha_{k}
))\notag\\
 =&S(S(\alpha_1)S(\alpha_2)\cdots S(\alpha_{l-1}) S(\beta_1) 
S(\beta_2) 
\cdots S(\beta_{r+1})) \, S(S(\alpha_{l+1}) \cdots S(\alpha_{k} )) \, n. \notag
\end{align}
Therefore, $(\alpha_1,\dots, \alpha_{l-1}, \beta_1, \dots, \beta_{r+1})$ and $(\alpha_{l+1},\dots,\alpha_k)$ are
both two-stack-sortable tuples of $r$-permutations.

Let us first consider the case  
when $l=1$. By (\ref{decomposition}), we have 
$$I=S(S(\beta_1) S(\beta_2) \cdots S(\beta_{r+1})) \, S(S(\alpha_{2})
\cdots S(\alpha_{k} )) \, n.$$ 
Suppose there are
$n-i$ letters (including $n$) in $\alpha_1$. Then the weight of
$(\alpha_2,\alpha_3, \dots, \alpha_{k})$ 
is a term of $g_i^{(k-1)}$. 
For $(\beta_1,\dots,
\beta_{r+1})$, if all $\beta_j$ are nonempty, then the weight of
$(\beta_1,\dots,\beta_{r+1})$ 
is a term of $g_{n-1-i}^{(r+1)}$. But by the
definition of the descents of $r$-permutations, there is a descent of each type in
$\alpha_1$ for every appearance of $n$. So the weight of $\alpha_1$ is $x_0x_1\cdots
x_r$ times the weight of $(\beta_1,\dots,\beta_{r+1})$. Therefore the sum of the
weights of all such $\alpha_1$ is $e_{r+1}(x)g_{n-1-i}^{(r+1)}$.
If only one of the 
$\beta_j$ is empty, then the weight of
$(\beta_1,\dots,\beta_{r+1})$ 
is a term of $g_{n-1-i}^{(r)}$. Again, by the definition of the
descents of $r$-permutations, there are $r-1$ descents of different types for $r-1$
 appearances of $n$. So, if $\beta_{j_0}$ is empty but $\beta_j$ is not empty for 
$j \ne j_0$, then the weight of $\alpha_1$ is $x_0 x_1 \cdots {\hat x_{j_0}} \cdots x_r$
times the weight of  $(\beta_1,\dots,\beta_{r+1})$. Therefore the sum of the weights of
all such
$\alpha_1$ is $e_r(x)g_{n-1-i}^{(r)}$. In general, if there are $m$ of the $\beta_j$
that are nonempty, then the sum of the weights of all such $\alpha_1$ is
$e_m(x)g_{n-1-i}^{(m)}$. 
Therefore, the sum of weights of all the $k$-tuple $r$-permutation 
$(\alpha_1,\alpha_2,\dots, \alpha_{k})$ in which $n$ is in $\alpha_1$ is
$$ \sum_{i=0}^{n-1}
g_i^{(k-1)}[g_{n-1-i}^{(0)}+g_{n-1-i}^{(1)}
e_1(x)+\cdots+g_{n-1-i}^{(r+1)}e_{r+1}(x)],$$

\noindent For the same reason,  when $n$ is in $\alpha_2$, we have 
$$I=S(S(\alpha_1)S(\beta_1) S(\beta_2) \cdots S(\beta_{r+1})) \,
S(S(\alpha_{3})
\cdots S(\alpha_{k} )) \, n.$$
Suppose there are $n-i$ letters in $(\alpha_1, \alpha_2)$. Then the sum of the weights
of all the $k$-tuple $r$-permutation 
$(\alpha_1,\alpha_2,\dots, \alpha_{k})$ in which $n$ is in $\alpha_2$ is 
$$ \sum_{i=0}^{n-1}
g_i^{(k-2)}[g_{n-1-i}^{(1)}+g_{n-1-i}^{(2)}
e_1(x)+\cdots+g_{n-1-i}^{(r+2)}e_{r+1}(x)],$$
\noindent and so on.
Thus we have the crucial recursive formula
\begin{equation}
g_n^{(k)}=\sum_{i=0}^{n-1} \sum_{l=0}^{k-1} \sum_{j=0}^{r+1}
g_{n-1-i}^{(l+j)} g_i^{(k-1-l)}e_j(x)
\end{equation}
with the initial conditions that $g_0^{(k)}=\delta_{0k}$ and $g_n^{(0)}=\delta_{n0}$.

Now  let
\begin{equation} \label{recursive1}
f_n^{(k)}=\sum_{j=0}^{r+1} g_n^{(k+j)} e_j(x).
\end{equation}
Then 
\begin{equation} \label{recursive2}
g_n^{(k)}=f_{n-1}^{(k-1)}+\sum_{i=0}^{n-1} f_{n-1-i}^{(k-2)} g_i^{(1)}+ \cdots
+\sum_{i=0}^{n-1} f_{n-1-i}^{(0)}g_i^{(k-1)}.
\end{equation}
Following the
analysis above one can easily see that $f_n^{(k)}$ is the sum of the weights
of 
$(k+1)$-tuple
$r$-permutations $(\alpha_1,\alpha_2, \dots ,\alpha_{k+1})$ of $[n+1]$ such 
that the $\alpha_i$
are all nonempty,   
$S(S(\alpha_1)S(\alpha_2)\cdots S(\alpha_{k+1}))=I$ and  the
largest letter only appears in $\alpha_{k+1}$. 

Also let
\begin{align} 
f_k=&\sum_{n=0}^{\infty}f_n^{(k)}, &F=\sum_{k=0}^{\infty}f_k z^k , \label{F} \\
g_k=&\sum_{n=0}^{\infty}g_n^{(k)}, &G=\sum_{k=0}^{\infty}g_k z^k . \label{G}
\end{align}

Then it follows that 
\begin{align}
g_0=&1, \notag \\
f_0=&g_1, \notag \\
f_k=&g_k+g_{k+1}e_1(x)+\cdots+g_{r+k+1}e_{r+1}(x), \label{FE1}\\
g_k=&f_{k-1}g_0+f_{k-2}g_1+\cdots+f_0g_{k-1}, \label{FE2}
\end{align}
and from (\ref{FE2}), we get
\begin{equation} \label{FE3}
G=1+zFG. 
\end{equation} 
Notice that if we let $r$ go to infinity, we can interpret
$k$-tuple $r$-permutations in terms of forests of 
decreasing trees and equations (\ref{FE1}) to (\ref{FE3}) are still well defined.
Therefore, we have the functional equations as follows.

\begin{Thm} \label{functional equations}
Let $F$ and $G$ be defined as in (\ref{F}) and (\ref{G}). Then $F$ and $G$ satisfy
the functional equations 
\begin{align}
f_k=&\sum_{j=0}^{\infty} g_{k+j}e_j(x) \text{ for $k \ge 0$}, \label{FE4}\\
 G=&1+zFG \label{FE5}
\end{align}
and $F$ and $G$ are uniquely determined as power series by these
functional equations.
\end{Thm}

\begin{Remark}
From (\ref{recursive2}), we see that $g_n^{(k)}$ are uniquely determined by
$f_i^{(j)}$ and
$g_i^{(j)}$ where $i<n$. From (\ref{recursive1}), we see that $f_n^{(k)}$ are
uniquely determined by $g_n^{(j)}$. Therefore, given $g_0^{(0)}=1$ and
$g_0^{(k)}=0$ when $k \ge 1$, all $g_n^{(k)}$ and $f_n^{(k)}$ are uniquely
determined by (\ref{recursive1}) and (\ref{recursive2}). Hence, $F$ and $G$ are
uniquely determined by the functional equations in Theorem \ref{functional equations}.
\end{Remark}
\begin{Remark}
Notice that  $g_1=f_0$ is the generating function of two-stack-sortable
$r$-permutations.
\end{Remark}
\begin{Remark}
Notice that in the ring ${\mathbb Q}[[z, \frac {x_0}{z}, \frac {x_1}{z}, \frac
{x_2}{z}, \dots]]$,
\begin{equation*}
G \cdot E\left(x,\frac{1}{z} \right) = \sum_{k=-\infty}^{\infty} 
\sum_{j=0}^{\infty} g_{k+j} e_j(x) z^k.
\end{equation*}
Therefore, by ({\ref{FE4}}), $F-G \cdot E\left(x,\frac{1}{z} \right)$ has only negative
powers of
$z$. This is the key to solving the functional equations (\ref{FE4}) and (\ref{FE5}).
\end{Remark}

\section{The solution to the functional equations}

Let
$y=(y_0,y_1,y_2,..)$, where
the $y_i$ are uniquely determined as formal power series in the $x_i$ by 
$$
{x_j={y_j(1+y_j)   \over
\prod_{i=0}^{\infty} (1+y_i)^2}={y_j(1+y_j) \over E^2(y)}}.$$ 
Let $t=zE^2(y)$ where $E(y)$ is defined as in (\ref{defE1}) and (\ref{E(x)}). Let
$\displaystyle {c(t)=\frac {1-\sqrt{1-4t}}{2t}}$ be the generating function of the
Catalan numbers. Then

\begin{Thm} \label{solution}
The solution to the functional equations in Theorem \ref{functional
equations} is:
\begin{align}
G=&\frac {c(t)E(y)}{E(y,c(t))}, \label{solution1}\\
F=&\frac{E(y)}{t}\left(E(y)-\frac{E(y,c(t))}{c(t)}\right). \label{solution2}
\end{align}
\end{Thm}

The proof of the theorem consists of the following lemmas.

\begin{Lem} \label{identity}
In the ring ${\mathbb
Q}[[z, \frac {x_0}{z}, \frac {x_1}{z}, \dots ]]$,
\begin{equation*} 
E\left(x,\frac{1}{z}\right)=E(y,c(t))E\left(y,\frac{1}{tc(t)}\right)
\end{equation*} 
\end{Lem}
\begin{proof}
First, we notice that $\displaystyle E\left(x,\frac {1}{z}\right) \in {\mathbb
Q}[[z, \frac {x_0}{z}, \frac {x_1}{z}, \dots ]]$ and  $\displaystyle
E\left(y,\frac{1}{tc(t)}\right) \in {\mathbb Q}[[t, \frac {y_0}{t}, \frac {y_1}{t},
\dots]]$. It is also easy to verify that $\displaystyle {\mathbb
Q}[[z, \frac {x_0}{z}, \frac {x_1}{z}, \dots ]]={\mathbb Q}[[t, \frac {y_0}{t}, \frac {y_1}{t},
\dots]]$. Therefore,
since $c(t)$ satisfies
$tc^2(t)-c(t)+1=0$, 
\begin{align*}
(1+c(t)y_i)\left(1+\frac{y_i}{tc(t)}\right)
=&1+c(t)y_i+\frac{y_i}{tc(t)}+\frac{y^2_i}{t} \\
=&1+\frac{y_i+y^2_i}{t}\\
=&1+\frac{x_i}{z}.
\end{align*}
Hence, 
\begin{align*}
E\left(x,\frac{1}{z}\right)=&\prod_{i=0}^{\infty}\left(1+\frac{x_i}{z} \right) \\
=&\prod_{i=0}^{\infty}\left( 1+\frac{y_i+y_i^2}{t}\right)\\
=&\prod_{i=0}^{\infty}(1+c(t)y_i)\left(1+\frac{y_i}{tc(t)}\right)\\
=&E(y,c(t))E\left(y,\frac{1}{tc(t)}\right).
\end{align*}
\end{proof}

\begin{Lem} \label{proveFE1}
Let  
\begin{align*}
\hat G=&\frac {c(t)E(y)}{E(y,c(t))},   \\
\hat F=&\frac{E(y)}{t}\left(E(y)-\frac{E(y,c(t))}{c(t)}\right).
\end{align*}
Then $\hat G$ and $\hat F$ satisfy $$\hat G=1+z \hat F \hat G.$$
\end{Lem}

\begin{proof}
\begin{align*}
1-z\hat F =& 1- z \cdot \frac{E(y)}{t}\left(E(y)-\frac{E(y,c(t))}{c(t)}\right) \\
=& 1- \frac {zE^2(y)}{t}\left( 1- \frac{E(y,c(t))}{c(t)E(y)}\right) \\
=& \frac{E(y,c(t))}{c(t)E(y)} \\
=& \frac{1}{\hat G}
\end{align*}
\end{proof}

\begin{Lem}  \label{negativeF} 
In the ring ${\mathbb Q}[[z, \frac {x_0}{z}, \frac
{x_1}{z}, \frac {x_2}{z}, \dots]]$, define 
$$\hat F^-=\hat G \cdot E\left(x, \frac{1}{z}
\right) - \hat F.$$ 
Then  
$\hat F^-$ has only negative powers of $z$.
\end{Lem}

\begin{proof}
{\allowdisplaybreaks
\begin{align*}
\hat F^-=&\hat G \cdot E\left(x, \frac{1}{z} \right) - \hat F \\
=&\frac {c(t)E(y)}{E(y,c(t))}\cdot E\left(x, \frac{1}{z}\right)- \frac{E(y)}{t}
\left(E(y)-\frac{E(y,c(t))}{c(t)}\right) \\
=&\frac {c(t)E(y)}{E(y,c(t))} \cdot E(y,c(t))E\left(y,\frac{1}{tc(t)}\right) - \frac{E(y)}{t}
\left(E(y)-\frac{E(y,c(t))}{c(t)}\right) \\
=& \frac{E(y)}{t}
\left[tc(t)E\left(y,\frac{1}{tc(t)}\right)+\frac{E(y,c(t))}{c(t)}-E(y)\right] \\
=& {E(y) \over t}
\sum_{n=0}^{\infty}\left[c(t)^{n-1}+t^{1-n}c(t)^{1-n}-1\right]e_n(y)\\
=& E(y) \sum_{n=0}^{\infty} \left[ \frac{1}{t^n} \left[ (tc(t))^{n-1}+c(t)^{1-n} \right]-\frac{1}{t} \right]
e_n(y).
\end{align*}
}
Since 
\begin{align}
(tc(t))^{n-1}+c(t)^{1-n}
=&\left( \frac {1-\sqrt{1-4t}} {2} \right)^{n-1} +
\left( \frac {1+\sqrt{1-4t}} {2} \right)^{n-1}  \\
=&\frac {1}{2^{n-2}} \sum_{k \text{ even}} \left( \sqrt{1-4t}\right)^k \binom{n-1}{k} \notag \\
=&\frac {1}{2^{n-2}} \sum_{k \text{ even}} (1-4t)^{k/2} \binom{n-1}{k} \notag
\end{align}
is a polynomial of degree less than or equal to $(n-1)/2$,
 $\hat F^-$ has only negative powers of $t$ and thus has only
negative powers of
$z$.
\end{proof}

\begin{Lem} \label{proveFE2}
Let $\hat G$ and $\hat F$ be defined as above. Define 
\begin{align*}
&\hat f_k= [z^k] \hat F, &\hat g_k= [z^k] \hat G.
\end{align*}
Then $\hat g_k$ and $\hat f_k$ satisfy the functional equation that 
\[
\hat f_k=\sum_{j=0}^{\infty}\hat g_{k+j}e_j(x) \text{ for $k \ge 0$}.
\]
\end{Lem}

\begin{proof}
From the definition of $\hat f_k$ and $\hat g_k$, we know that  
\begin{align*}
\hat F=& \sum_{k=0}^{\infty} \hat f_k z^k, \\
\hat G=& \sum_{k=0}^{\infty} \hat g_k z^k.
\end{align*}
Then 
$$ [z^k] \left\{ \hat G \cdot E\left( x, \frac {1}{z} \right) - \hat F
\right\} =\sum_{j=0}^{\infty}\hat g_{k+j}e_j(x)-\hat f_k . $$
By Lemma \ref{negativeF}, this sum is zero when $k \ge 0$.

The lemma is proved.
\end{proof}

Now Theorem \ref{solution} follows from Lemma \ref{proveFE1} and Lemma \ref{proveFE2}.

\section{The number of two-stack-sortable $r$-permutations with descents}

 We know that the generating function for two-stack-sortable
$r$-permutations is 
{\allowdisplaybreaks
\begin{align*}
g_1=&f_0 \\
=&[z^0]F \\
=&[t^0] \frac{E(y)}{t}\left(E(y)-\frac{E(y,c(t))}{c(t)}\right) \\
=&[t^0] \left\{E(y)  \sum_{n=0}^{\infty} \frac{1}{t} (1-c^{n-1}(t))e_n(y)
\right\}\\
=&[t^0]\left\{ E(y)\sum_{n=0}^{\infty} \sum_{i=0}^{\infty} \frac{1-n}{i+1}
\binom{n+2i}{i}e_n(y) t^i\right\} \\
 =& E(y) \sum_{n=0}^{\infty} (1-n)e_n(y)
\end{align*}
}
where 
\begin{equation*}
x_i=\frac {y_i(1+y_i)}{E^2(y)}.
\end{equation*} 
We can now evaluate the number of two-stack-sortable $r$-permutations with the help
of multivariable Lagrange inversion. 

\begin{Thm} \label{numberwithdescents}
The number of two-stack-sortable $r$-permutations with $k_i$ descents of
type
$i$ is 
\begin{equation}
\frac{1}{n^2}\prod_{i=0}^r \frac{n}{n-k_i} \binom {2n-1-k_i} {k_i},
\end{equation}
where $n=1+k_0+k_1+\cdots+k_r$.
\end{Thm}

First, we introduce the following identity:

\begin{Lem}
$$\sum_{i=0}^{r+1} (1-i)e_i(y)=E(y) \left(1-\sum_{i=0}^r \frac{y_i}{1+y_i} 
\right).$$ 
\end{Lem}
\begin{proof}
Since $e_i(y)$ is a homogeneous polynomial in $y_0,y_1,\dots,y_r$ of degree $i$, 
$$ i e_i(y)=\sum_{j=0}^r \frac {y_j \partial e_i(y)} {\partial y_j}.$$
Also notice that 
$$ \frac {\partial E(y)} {\partial y_j} =\frac {E(y)}{1+y_j}. $$
Thus 
{\allowdisplaybreaks
\begin{align*}
\sum_{i=0}^{r+1} (1-i)e_i(y)=& E(y)-\sum_{i=0}^{r+1} ie_i(y)\\
=&E(y)-\sum_{i=0}^{r+1}\sum_{j=0}^r \frac {y_j \partial e_i(y)}
{\partial y_j}\\
=&E(y)-\sum_{j=0}^r y_j \sum_{i=0}^{r+1}\frac { \partial e_i(y)}
{\partial y_j} \\
=& E(y) -\sum_{j=0}^r y_j \frac {E(y)}{1+y_j} \\
=&E(y) \left(1-\sum_{j=0}^r \frac{y_j}{1+y_j}  \right).
\end{align*}
}
\end{proof}

Therefore, 
\begin{align} \label{g1}
g_1=& E(y) \sum_{n=0}^{\infty} (1-n)e_n(y) \\
=&E^2(y) \left(1-\sum_{j=0}^r \frac{y_j}{1+y_j}  \right). \nonumber
\end{align}

If for any $y=(y_0,y_1,\dots,y_r)$, we define 
\begin{align*}
A(y,u)=& A(y_0,\dots,y_r;u) \\
=&1-\sum_{i=0} \frac{y_iu}{1+y_iu},
\end{align*}
and let $$A(y)=A(y,1),$$
then $g_1=E^2(y)A(y)$.

The proof of Theorem \ref{numberwithdescents} now consists of following lemmas.

Notice that for any formal power series $f(x)= \sum_k a_k x^k \in 
{\mathbb Q} [[ x_0, x_1, \dots, x_r ]]$ 
(where $k=(k_0,\dots,k_r)$ and $x=(x_0,\dots,x_r)$), if we define operators $D_x$, $I_x$:
  ${\mathbb Q} [[ x_0, x_1, \dots, x_r ]] \to
{\mathbb Q} [[ x_0, x_1, \dots, x_r]]$ by 
\begin{align*}
D_x (f(x))=& \sum_{i=0}^r x_i\frac{\partial f(x)}{\partial x_i} , \\
I_x (f(x))=& f(x),
\end{align*}
then  $(D_x+I_x)(f(x))=\sum_k n a_kx^k$, where $n=k_0+\cdots+k_r+1$.
Therefore, to prove
$$[x^k]g_1=\frac {1}{n^2} \prod_{i=0}^r \frac{n}{n-k_i} \binom {2n-1-k_i} {k_i},$$
we only need to prove that 
$$[x^k](D_x+I_x)^2g_1=\prod_{i=0}^r \frac{n}{n-k_i} \binom {2n-1-k_i} {k_i}.$$

Let $u_i=\sum_{j=0}^r \dfrac{x_j \partial y_i}{\partial x_j}$. Then

\begin{Lem}
$$u_i=\frac {y_i(1+y_i)}{1+2y_i} \frac{1}{A(y,2)}$$
\end{Lem}
\begin{proof}
Take the logarithm of both sides of $x_i=\dfrac{y_i(1+y_i)}{E^2(y)}$. We get
\begin{equation} \label{log}
\ln x_i= \ln y_i + \ln (1+y_i) - 2 \ln E(y).
\end{equation}
Then by differentiating both sides of \ref{log} with respect to $x_j$, we get
\begin{align}
\frac{1}{x_i}=&\frac{1+2y_i}{y_i(1+y_i)} \frac {\partial y_i}{\partial x_i}
-\frac{2}{E(y)} \frac {\partial E(y)}{\partial x_i}, & \text{ if   $j=i$}, \label{d1}\\
0=& \frac{1+2y_i}{y_i(1+y_i)} \frac {\partial y_i}{\partial x_j}
-\frac{2}{E(y)} \frac {\partial E(y)}{\partial x_j}, & \text{ if   $j \ne i$}.
\label{d2}
\end{align}
Now multiplying $x_j$ by both (\ref{d1}) and (\ref{d2}) and summing on $j$, we get 
\begin{align}
1=&\frac{1+2y_i}{y_i(1+y_i)} \sum_{j=0}^r \frac {x_j\partial y_i}{\partial x_j}
-\frac{2}{E(y)} \sum_{j=0}^r\frac {x_j\partial E(y)}{\partial x_j} \nonumber\\
=&\frac{1+2y_i}{y_i(1+y_i)} u_i -\frac{2}{E(y)} D_x(E(y)). \label{fun1}
\end{align}

On the other hand,
\begin{align}
D_x(E(y))=& \sum_{j=0}^r x_j \frac{\partial E(y)}{\partial x_j} \nonumber \\
=&\sum_{j=0}^r x_j \sum_{i=0}^r \frac{\partial  E(y)}{\partial y_i}
\frac{\partial y_i}{\partial x_j} \nonumber \\
=&E(y) \sum_{i=0}^r \frac{1}{1+y_i} \sum_{j=0}^r x_j 
\frac{\partial y_i}{\partial x_j} \nonumber \\
=&E(y) \sum_{i=0}^r \frac{1}{1+y_i} u_i \label{fun2}
\end{align}
By solving (\ref{fun1}) and (\ref{fun2}) for $u_i$, we get 
$$u_i=\frac {y_i(1+y_i)}{1+2y_i} \frac{1}{A(y,2)}.$$
The lemma is proved.
\end{proof}

Using  the lemma above, the proofs of the following lemmas are just a matter of simple
computations.
\begin{Lem}\label{identity1}
\[
D_x(E^2(y))=\frac{E^2(y)}{A(y,2)}-E^2(y).
\]
\end{Lem}
\begin{proof}
{\allowdisplaybreaks
\begin{align*}
D_x(E^2(y))
=&\sum_{j=0}^r x_j \frac{\partial E^2(y)}{\partial x_j}  \\
=&\sum_{j=0}^r x_j\sum_{i=0}^r 
\frac {\partial E^2(y)}{\partial y_i} \frac {\partial y_i}{\partial x_j}  \\
=&\sum_{i=0}^r\frac {\partial E^2(y)}{\partial y_i}  u_i  \\
=&\sum_{i=0}^r \frac{2E^2(y)}{1+y_i} \frac {y_i(1+y_i)}{1+2y_i} \frac{1}{A(y,2)}\\
=&\frac{E^2(y)}{A(y,2)} \sum_{i=0}^r \frac{2y_i}{1+2y_i}  \\
=&\frac{E^2(y)}{A(y,2)}-E^2(y).
\end{align*}
}
The lemma is proved.
\end{proof}
\begin{Lem} \label{identity2}
\[
D_x(A(y))=1-\frac {A(y)}{A(y,2)}.
\]
\end{Lem}
\begin{proof}
{\allowdisplaybreaks
\begin{align*}
D_x(A(y))
=&\sum_{j=0}^r x_j \frac{\partial A(y)}{\partial x_j}  \\
=&\sum_{j=0}^r x_j\sum_{i=0}^r 
\frac {\partial A(y)}{\partial y_i} \frac {\partial y_i}{\partial x_j}  \\
=&\sum_{i=0}^r\frac {\partial A(y)}{\partial y_i}  u_i  \\
=& \sum_{i=0}^r -\frac{1}{(1+y_i)^2} \frac {y_i(1+y_i)}{1+2y_i} \frac{1}{A(y,2)} \\
=&\frac{1}{A(y,2)} \sum_{i=0}^r - \frac{y_i}{(1+y_i)(1+2y_i)}  \\
=&\frac{1}{A(y,2)}  \sum_{i=0}^r \left( \frac{y_i}{1+y_i} -\frac{2y_i}{1+2y_i} \right)\\
=&\frac{1}{A(y,2)} \left(A(y,2) -A(y)\right)  \\
=&1-\frac {A(y)}{A(y,2)}.
\end{align*}
}
\end{proof}
\begin{Lem} \label{identity3}
\[(D_x+I_x)(g_1)=E^2(y)\]
\end{Lem}
\begin{proof}
\begin{align*}
(D_x+I_x)(g_1)
=&(D_x+I_x) (E^2(y)A(y))\\
=&E^2(y) D_x (A(y))+A(y) D_x(E^2(y))+E^2(y)A(y)\\
=&E^2(y) \left(1-\frac {A(y)}{A(y,2)} \right) + 
A(y) \left( \frac{E^2(y)}{A(y,2)} -E^2(y) \right) +E^2(y)A(y)\\
=&E^2(y).
\end{align*}
\end{proof}

Now from Lemma \ref{identity1} to Lemma \ref{identity3}, we see that 
\[
\frac {E^2(y)}{A(y,2)} = (D_x+I_x)^2(g_1).
\]

\begin{Lem} \label{identity4}
\[
[x^k]\frac {E^2(y)}{A(y,2)}
=  \prod_{i=0}^r \frac {n}{n-k_i} \binom{2n-1-k_i}{k_i} 
\] where $n = k_0+\cdots+k_r+1$.
\end{Lem}
\begin{proof}
Using multivariable Lagrange inversion (Lemma \ref {lagrange}), we set
$$\phi_i=\frac{E^2(y)}{1+y_i}.$$ Thus $y_i=x_i \phi_i(y)$.
Using Lemma \ref{determinant} to evaluate the determinant, we have
{\allowdisplaybreaks
\begin{align*}
[x^k]\frac {E^2(y)}{A(y,2)}
=& [y^k] \left\{ \frac {E^2(y)}{A(y,2)} \phi^k(y)
\left| \delta_{ij} -\frac{y_j}{\phi_i(y)} \frac{\partial
\phi_i(y)}{\partial y_j} \right| \right\} \\
=&  [y^k] \left\{ \frac {E^2(y)}{A(y,2)} \phi^k(y) 
\frac {\left| \delta_{ij}+2(\delta_{ij}-1)y_j \right|} {E(y)} \right\} \\
=&[y^k] \left\{ \frac {E^2(y)}{A(y,2)} \phi^k(y) 
\frac {E(y,2)A(y,2)}{E(y)} \right\} \\
=& [y^k] \left\{  \frac
{E^{2n-1}(y)} {\prod_i (1+y_i)^{k_i}}  E(y,2) \right\}\\
=& \prod_i \left[ \binom
{2n-1-k_i}{k_i}+2\binom{2n-1-k_i}{k_i-1} \right] \\
=&\prod_i \frac{n}{n-k_i} \binom {2n-1-k_i} {k_i}. 
\end{align*}
}
\end{proof}

Therefore, Theorem \ref{numberwithdescents} is proved.

\vs 4

\begin{Remark} 
When $r=1$, we have $k_0+k_1=n-1$. Then the number is 
\begin{align*}
\frac {1}{n^2} \frac {n}{n-k_0} \binom{2n-1-k_0}{k_0} 
&\frac {n}{n-k_1} \binom{2n-1-k_1}{k_1} \nonumber\\
=& \frac {(n+k_1)! \, (2n-1-k_1)!}
{(k_1+1)! \, (2k_1+1)!\, (n-k_1)! \, (2n-1-2k_1)!},
\end{align*}
which was found and proved by B{\'o}na [\ref{Bona2}], Bousquet-M{\'e}lou 
[\ref{melou3}], and Travis [\ref{Travis}].
\end{Remark}

\section{The number of two-stack-sortable $r$-permutations}

The following theorem is a special case of Theorem
\ref{numberwithdescents} when  we do not keep track of descents.
 
\begin{Thm} 
The number of two-stack-sortable $r$-permutations of $[n]$ is:
$$
2(r+1) \frac {((2r+1)n)!} {n! \, (2rn+2)!}.
$$
\end{Thm}

\noindent {\it First proof:}

We know that the generating function for two-stack-sortable $r$-permutations
is given by (\ref{g1}).
Setting $y_i=y$ for $i=0,\dots,r$, which is equivalent to setting $x_i=x$ for
$i=0,\dots,r$, we get that 
$$x=\frac {y} {(1+y)^{2r+1}}$$ and
\begin{align}
g_1=&(1+y)^{2r+2} \left( 1- (r+1)\frac {y}{1+y} \right) \notag \\
=& (r+1) (1+y)^{2r+1} -r (1+y)^{2r+2} \notag \\
=& (1+y)^{2r+1}(1-ry) \notag \\
=& \frac {y(1-ry)}{x}. \label{g1expression}
\end{align}

Using Lagrange inversion [\ref{Goulden}, p.~ 17; \ref{Ira}] , we set $\phi(y)=(1+y)^{2r+1}$. Then
$y$ satisfies the condition that $y=x \phi (y)$.
Therefore, 
\begin{align*}
[x^{k}]g_1=&[x^{k+1}] xg_1 \\
=& \frac{1}{k+1} [y^{k}] \left\{ \frac{d (xg_1)}{d y} \phi^{k+1}(y)
\right\} \\
=&\frac{1}{k+1} [y^{k}] \left\{ (1-2ry) 
(1+y)^{(2r+1)(k+1)} \right\} \\
=& \frac{1}{k+1} \left\{  \binom {(2r+1)(k+1)}{k} - 2r
\binom {(2r+1)(k+1)}{k-1} \right \} \\
=&2(r+1) \frac {((2r+1)(k+1))!} {(k+1)! \, (2r(k+1)+2)!}.
\end{align*}
Setting $n=k+1$, then the theorem  is proved.

\vs  4

\noindent {\it Second proof:}

From Theorem \ref{numberwithdescents}, we know that the total number of
two-stack-sortable permutations is 
\[
\frac{1}{n^2} \sum_{k_0+\cdots+k_r=n-1} \prod_{i=0}^r \frac {n}{n-k_i} 
\binom{2n-1-k_i}{k_i}.
\]
Also we know that  
\begin{equation} \label{c(t)}
c^n(t)=\sum_{i=0}^{\infty} \frac{n}{2i+n} \binom{2i+n}{i}t^i 
\end{equation}
for any integer
$n$  [\ref{john}, p.~ 154].
Therefore,
\begin{align*}
[t^{k_i}]c^{-2n}(-t)=& [t^{k_i}] \sum_{j=0}^{\infty} \frac {2n}{2n-2j} \binom
{2j-2n}{j} (-t)^j \\
=& \frac  {2n}{2n-2k_i} (-1)^{k_i} \binom {2k_i-2n}{k_i} \\
=& \frac{n}{n-k_i} \binom{2n-1-k_i}{k_i}.
\end{align*}
Hence, 
\begin{align*}
\frac{1}{n^2} \sum_{k_0+\cdots+k_r=n-1} \prod_{i=0}^r \frac {n}{n-k_i} 
&\binom{2n-1-k_i}{k_i} \\
=& \frac{1}{n^2} [t^{n-1}] \left( c^{2n}(-t) \right) ^{r+1} \\
=& \frac{1}{n^2} [t^{n-1}] \sum_{j=0}^{\infty} 
    \frac {2n(r+1)} {2n(r+1)-2j} \binom {2j- 2n(r+1)}{j} (-t)^j \\
=& \frac{1}{n^2} \frac {2n(r+1)}{2n(r+1)-2(n-1)} (-1)^{n-1} 
				\binom {2(n-1)-2n(r+1)}{n-1} \\
=& \frac {r+1}{n(nr+1)} \binom {(2r+1)n}{n-1} \\
=& 2(r+1) \frac {((2r+1)n)! } {n! \, (2rn+2)! }.
\end{align*}
The theorem is proved.
\qed

\vs 5
\begin{Remark} When $r=1$, this number is  
\[ 
4 \frac {(3n)!}{ n! \, (2n+2)!}=\frac {2(3n)!}{(n+1)! \, (2n+1)!},
\]
which was conjectured by West [\ref{West2}, \ref{West}] and first proved by
Zeilberger [\ref{Zeil}].
\end{Remark}

\chapter{Parallel results}

\section{Another approach} \label{another approach}
In the previous chapter, we counted the number of two-stack-sortable $r$-permutations
with all types of descents. Here we approach a special case of the problem from another
angle. We now count the number of two-stack-sortable $r$-permutations with only 
two basic parameters, the number of different letters and the number of components. 
 We use the indeterminates $x$ and $\bar z$ to count these two parameters respectively. The
difference is, that {\em any component is allowed to be empty}. 

Therefore, using the same decomposition method as in Chapter 2, if we let $p_n^{(k)}$ be the
weight of all 
$k$-tuple $r$-permutations $(\alpha_1, \dots, \alpha_k)$ of $[n]$ where each of $\alpha_i$
$(i=1,\dots,k)$ is {\em allowed to be empty}, and let  
$q_n^{(k)}=p_n^{(r+k+1)}$,
\begin{align} 
p_k=&\sum_{n=0}^{\infty}p_n^{(k)}x^n, &P=\sum_{k=0}^{\infty}p_k {\bar z}^k , \label{P} \\
q_k=&\sum_{n=0}^{\infty}q_n^{(k)}x^n, &Q=\sum_{k=0}^{\infty}q_k {\bar z}^k , \label{Q}
\end{align}
then analogous to Theorem \ref{functional equations}, we have
that the formal power series $P$ and $Q$ satisfy the functional equations
\begin{align}
q_k=&p_{r+k+1}, \label{FE21}\\
P=&\frac{1}{1-\bar z}+x{\bar z}PQ, \label{FE22}
\end{align}
Since $$Q=\dfrac {1}{ {\bar z}^{r+1}}\left(P-\sum_{i=0}^r p_i {\bar z}^i\right), $$ we have that
\begin{Thm}
The generating function $P$ of two-stack-sortable $k$-tuple $r$-permutations satisfies 
\begin{equation} \label{FE23}
P = \frac {1}{1-\bar z} +\frac {xP}{{\bar z}^r} \left(P-\sum_{i=0}^r p_i {\bar z}^i\right)
\end{equation}
and $P$ is uniquely determined as a power series by this functional equation.
\end{Thm}

Now, compared to the functional equations (\ref{FE4}) and (\ref{FE5}), this one is structurally
simpler. We will also see that (\ref{FE23}) can be easily derived from the functional
equations (\ref{FE4}) and (\ref{FE5}).

If we set $x_i = x$ in functional equations
(\ref{FE4}) and (\ref{FE5}), then they become
\begin{align}
f_k=&\sum_{j=0}^{r+1} \binom{r+1}{j} g_{k+j} x^j \text{ for $k \ge 0$}, \label{FE6}\\
 G=&1+zFG. \label{FE7}
\end{align}

Notice that $p_k$ is the weight of of all 
$k$-tuple  $r$-permutations  where {\em any component is  allowed to be empty} and $g_k$ is the 
weight of of all 
$k$-tuple $r$-permutations  where {\em no component is  allowed to be empty}. Therefore,
\begin{equation} \label{pkgk}
p_k=\sum_{i=0}^k \binom{k}{i} g_i x^i ,
\end{equation}
which is equivalent to 
$$ P(x, \bar z)=\frac{1}{1-\bar z}G\left(x, \frac{x \bar z}{1-\bar z} \right). $$
 
Setting $x_i=x$, which is equivalent to setting $y_i=y$, in the
solution of functional equations  (Theorem \ref{solution}) gives the solution to functional
equation (\ref{FE6}) and (\ref{FE7}), 
$$
G=\frac {c(t)(1+y)^{r+1}} { (1+ yc(t))^{r+1}},
$$
where $y$ is uniquely determined as formal power series in $x$ by 
$$ 
x=\frac {y}{(1+y)^{2r+1}}
$$
and $t=z(1+y)^{2(r+1)}$. Therefore,
\begin{Thm} \label {solution2}
The solution to function equation (\ref{FE23}) is:
$$
P=\left(1+\frac{t}{y(1+y)}\right) \frac {c(t)(1+y)^{r+1}} { (1+ yc(t))^{r+1}},
$$
where $y$ is uniquely determined as a formal power series in $x$ by 
$$ 
x=\frac {y}{(1+y)^{2r+1}}
$$
and $t=y(1+y)\dfrac{\bar z}{1-\bar z}$.
\end{Thm}

\begin{Cor} \label{pkpoly}
$p_k$ is a polynomial in $y$ with degree no greater than $2k$.
\end{Cor}
\begin{proof}
{\allowdisplaybreaks
\begin{align*}
x^kg_k=&x^k [z^k]G \\
=&\left[\left(\frac{z}{x}\right)^k \right]G \\
=&\left[\left(\frac{z}{x}\right)^k \right] \left\{\frac {c(t)(1+y)^{r+1}} { (1+
yc(t))^{r+1}} \right\} \\ 
=&[z^k] \left\{\frac {c(xt)(1+y)^{r+1}} { (1+
yc(xt))^{r+1}} \right\} \\ 
=&[{z}^k] \left\{(1+y)^{r+1} 
\sum_{i=0}^{\infty}(-1)^i \binom {r+i}{i}y^i (c(xt))^{i+1} \right\} \\
=& [{z}^k] \left\{(1+y)^{r+1} 
\sum_{i=0}^{\infty}(-1)^i \binom {r+i}{i}y^i
\sum_{j=0}^{\infty} \frac {i+1}{2j+i+1} \binom {2j+i+1}{j}(xt)^j \right\} \\
=&[{z}^k] \left\{(1+y)^{r+1} 
\sum_{i=0}^{\infty}(-1)^i \binom {r+i}{i}y^i
\sum_{j=0}^{\infty} \frac {i+1}{2j+i+1} \binom {2j+i+1}{j}
y^j(1+y)^j z^j\right\} \\
=&(1+y)^{r+k+1} y^k\sum_{i=0}^{\infty}(-1)^i \binom {r+i}{i}y^i
\frac {i+1}{2k+i+1} \binom {2k+i+1}{k}.
\end{align*}
}
The sum is a hypergeometric series that can be transformed by Euler's transformation [\ref{bailey}]. Thus,
$$x^kg_k=y^k \sum_{i=0}^k (-1)^i\frac {(2i+1)r-(k-i)}{(r-k)(k+i+1)}
\binom{r-k}{i} \binom {2k}{k-i} y^i.
$$
Therefore, $x^kg_k$ is a polynomial in $y$ with degree no greater than $2k$.
Thus $p_k$ is a polynomial in $y$ with degree no greater than $2k$ because 
$p_k=\sum_{i=0}^k \binom{k}{i}  x^ig_i$ (Eq.~ (\ref{pkgk})).
\end{proof}

In particular, 
\begin{Cor} \label{p_1} 
$p_1$ is the generating function of two-stack-sortable $r$-permutations, and 
$$p_1=1+y-ry^2.$$
\end{Cor}
\begin{proof}
Since $p_1=1+xg_1$ (Eq. (\ref{pkgk})) and $g_1=\dfrac {y(1-ry)}{x}$ (Eq. (\ref{g1expression})),
\begin{align*}
p_1=& 1 + x g_1 \\
=& 1+ y -ry^2.
\end{align*}
\end{proof}

\section{Connection to Zeilberger's  functional equation}

Zeilberger first proved West's conjecture that the number of two-stack-sortable
permutations of length $n$ is $ \dfrac {2(3n)!}{(n+1)! \, (2n+1)!}$ [\ref{Zeil}]. He
used a factorization similar to ours and the functional equation he got is equivalent
to our functional equation (\ref{FE23}) in the case of $r=1$. 

In Zeilberger's paper [\ref{Zeil}], he defined $i(\pi)$ (where $\pi$ is any permutation
of 
$\{1,2,\dots,n\}$) to be the largest integer $i$ such that the subsequence of the 
 `big $i$': $\{n-i+1, \dots, n-1, n\}$ are in decreasing order, defined $W^{(i)}$ to be
the set of all permutations (of any length) $\pi$ such that $i(\pi)=i$, and let 
$W^{\ge i}$ to be the set of all permutations $\pi$ such that $i(\pi) \ge i$. Also he 
defined 
$W^{(i)}(x)$ to be the formal power series that equals to the sum of all the weights of 
 elements of $W^{(i)}$, and $W^{\ge i}(x)$ to
be the formal power series that equals to the sum of all the weights  of
elements of $W^{\ge i}$, and he defined 
\begin{align*}
\Phi(x,t):=&\sum_{i=0}^{\infty}W^{(i)}(x)t^i, \\
{\bar \Phi(x,t)}:=&\sum_{i=0}^{\infty}W^{\ge i}(x)t^i.
\end{align*}
Then, he got 
\begin{equation} \label{zeil's}
\Phi(x,t)=\frac{1}{1-xt}+\frac{xt(\Phi(x,1)-t\Phi(x,t))(\Phi(x,1)-\Phi(x,t))}
{(1-t)^2}.
\end{equation}

Now noticing that $\Phi(x,1)={\bar \Phi(x,0)}$ and 
\[
{\bar \Phi}(x,t)=\frac{\Phi(x,1)-t\Phi(x,t)}{1-t}.
\]
Equation (\ref{zeil's}) is equivalent to
\[
{\bar \Phi}(x,t)=\frac{1}{1-xt}+
\frac {(1+xt {\bar \Phi}(x,t)) ({\bar \Phi}(x,t) -{\bar \Phi}(x,0) )} {t}.
\]
Now it is easy to check that $P(x,\bar z)=1+{\bar z} {\bar \Phi}\left(x,\dfrac{\bar z}{x}\right)$. 
Therefore, the functional equation (\ref{zeil's}) is equivalent to (\ref{FE23}) in the
case of $r=1$.

The combinatorial connection between the set $W^{\ge i}$ and the set of $k$-tuple
permutations was clearly stated in Zeilberger's paper [\ref{Zeil}]. For a typical
member $\pi$ of $W^{\ge i}$, if its length is $n$, then it has the form 
$$
\pi=\alpha_0 n \alpha_1 (n-1) \cdots (n-i+1) \alpha_i,
$$
where $\alpha_0,\dots,\alpha_i$ are (possibly empty) permutations of disjoint smaller
sets, the union of whose underlying sets is $\{1,2,\dots,n-i\}$. Now, by iterating the
definition of the stack sorting operation $S$,
$$
S(\pi)=S(\alpha_0) S(\alpha_1) \cdots S(\alpha_i) (n-i+1)(n-i+2) \cdots n,
$$
so that,
$$
S^2(\pi)=S(S(\alpha_0) S(\alpha_1) \cdots S(\alpha_i)) (n-i+1)(n-i+2) \cdots n.
$$
It follows that there is a 1-1 correspondence between the elements of $W^{\ge i}$ and 
$(i+1)$-tuple permutations $\alpha_0,\alpha_1,\dots, \alpha_i$, such that $S(S(\alpha_0)
 \cdots S(\alpha_i))=I$, and the underlying sets of the $\alpha$'s are disjoint and
their union is $\{1,2,\dots,n-i\}$.

\chapter{Further results}
\section{A generalization of the stack-sorting operation on $r$-permutations}
Here we introduce a more general form of the stack-sorting operation on
$r$-permutations.

\begin{Def}
For a given positive number $r$, let $\lambda=(\lambda_0,\lambda_1,\dots,
\lambda_l)$ satisfy the condition that 
$0 < \lambda_0 < \lambda_1 < \cdots < \lambda_l=r$ . 
Let $\pi$ be an
$r$-permutation of 
$[n]$. Then we can write $\pi$ as $\pi=\alpha_0 n \alpha_1 n \cdots n
\alpha_{r}$. 
The generalized stack-sorting operation
$S_{\lambda}$ on $r$-permutations is defined by 
$$S_{\lambda}(\pi)=S_{\lambda}(\alpha_0)\cdots S_{\lambda}(\alpha_{\lambda_0}) n
S_{\lambda}(\alpha_{\lambda_0+1})\cdots S_{\lambda}(\alpha_{\lambda_1}) n
\cdots n S_{\lambda}(\alpha_{\lambda_{l-1}+1})\cdots
S_{\lambda}(\alpha_{\lambda_l}) n.
$$
Notice that $S_{\lambda}(\pi)$ is an $(l+1)$-permutation and $S_{\lambda}=S$
when $l=0$, where $S$ is the ordinary stack
sorting operation on
$r$-permutations (Definition \ref{stacksortingoperation}).
\end{Def}

It is clear that it is not interesting to consider  stack-sortable
$r$-permutation under this definition. Also, in order to consider
two-stack-sortable $r$-permutations under this definition, we have to consider
$S(S_{\lambda}(\pi))$ instead of $S_{\lambda}^2(\pi)$.

To enumerate the number of two-stack-sortable $r$-permutations with descents under this
definition of stack-sorting operation, we still weight any type $i$ descent by $x_i$.
Furthermore, we denote $x^{(0)}=(x_0,x_1,\dots,x_{\lambda_0})$ and 
$x^{(i)}=(x_{\lambda_{i-1}+1},\dots, x_{\lambda_i})$ for
$1 \le i \le l$. Also, we use $u$ to keep track of the number of different letters.

Now let $g_n^{(k)}$ be the sum of the weights of two-stack-sortable $k$-tuple
$r$-permutations $(\alpha_1,\alpha_2,\dots, \alpha_k)$ such that every $\alpha_i$ is
nonempty and let 
\begin{align*}
f_{n,0}^{(k)}=&\sum_{j=0}^{\lambda_0+1} g_n^{(k+j)} e_j (x^{(0)} ), \\
f_{n,i}^{(k)}=&\sum_{j=0}^{\lambda_i-\lambda_{i-1}} g_n^{(k+j)} e_j(x^{(i)})
\, \text{   for  }  \, 1 \le i \le l.
\end{align*}

Also let 
\begin{align} \label{defofFG} 
f_{k,i}=&\sum_{n=0}^{\infty}f_{n,i}^{(k)}u^n, &&F_i=\sum_{k=0}^{\infty}f_{k,i} z^k ,\\
g_k=&\sum_{n=0}^{\infty}g_n^{(k)}u^n, &&G=\sum_{k=0}^{\infty}g_k z^k \notag. 
\end{align}

Then we have the following functional equations by the same reasoning as in section \ref{mainfuneq}:
\begin{Lem}
Let $G$ and $f_{k,i}$ be defined as above, then
\begin{align}
G=&1+z u G F_0 \prod_{i=1}^l f_{0,i}, \label{FE41}\\
f_{k,0}=&\sum_{j=0}^{\lambda_0+1} g_{k+j}e_j(x^{(0)}),\label{FE42}\\
f_{k,i}=&\sum_{j=0}^{\lambda_i-\lambda_{i-1}} g_{k+j} e_j(x^{(i)}) \, 
\mbox{  for  } 1 \le i \le l. \label{FE43}
\end{align}
\end{Lem}
If we let $A=\prod_{i=1}^l f_{0,i}$, we can make the substitution $x=u A $ so that the functional
equations  (\ref{FE41}) to (\ref{FE43}) become 
\begin{align}
G=&1+z x G F_0 , \\
f_{k,0}=&\sum_{j=0}^{\lambda_0+1} g_{k+j}e_j(x^{(0)}). 
\end{align}

By Theorem \ref{solution}, the solution to these functional equations is:
\begin{align*}
G=&\frac {c(t)E(y)}{E(y,c(t))}, \\
F_0=&\frac{E(y)}{t}\left(E(y)-\frac{E(y,c(t))}{c(t)}\right),
\end{align*}
where $y=(y_0,\dots,y_{\lambda_0})$ is uniquely determined by 
$xx_j=\dfrac {y_i (1+y_i)}{E^2(y)}$ for $0 \le j \le \lambda_0$, and $t=zxE^2(y)$.
From the solution, we can express $g_i$ in terms of $x$ and $y$, i.e., $g_i=g_i(x,y)$. Since $x=uA$, $g_i$ can
be expressed in terms of $u,A,y$, i.e., $g_i=g_i(u,A,y)$. By equation (\ref{FE43}), 
\begin{align*}
A=&\prod_{i=1}^l f_{0,i} \\
=&\prod_{i=1}^l\sum_{j=0}^{\lambda_i-\lambda_{i-1}} g_{j}(u,A,y) e_j(x^{(i)}).
\end{align*}
Now we can solve for $A$ from the equation above so that $A$ is expressed in terms of $u,y,x^{(i)}$, 
$1 \le i \le l$, i.e., $A=A(u;y;x^{(1)};\dots;x^{(l)})$.
Therefore, given $g_i=g_i(u,A,y)$ and $x_j=\dfrac {y_i
(1+y_i)}{xE^2(y)}$ for $0 \le j \le \lambda_0$, we can solve for $g_i$ (for example, using Lagrange inversion)
so that
$g_i=g_i(u;x^{(0)};x^{(1)};\dots;x^{(l)})$.

Therefore,
\begin{Thm} \label{solution3}
If $G$ and $F_i$ are defined in (\ref{defofFG}), then $G$ and $F_i$ satisfy the following functional
equations:  
\begin{align*}
G=&\frac {c(t)E(y)}{E(y,c(t))}, \\
F_0=&\frac{E(y)}{t}\left(E(y)-\frac{E(y,c(t))}{c(t)}\right),
\end{align*}
where $y=(y_0,\dots,y_{\lambda_0})$ is uniquely determined by 
$xx_i=\dfrac {y_i (1+y_i)}{E^2(y)}$ for $0 \le i \le \lambda_0$, and $t=zxE^2(y)$.
In particular, the $g_i$ are algebraic.
\end{Thm}
\begin{Remark} \label{g1ex}
It is easy to get $g_0=1$. Also, similar to Equation (\ref{g1}), we can get 
\[
g_1=x E(y) \sum_{n=0}^{\infty} (1-n)e_n(y)
\]
in this case.
\end{Remark}

For example, the simplest nontrivial case is when $l=1$ and $\lambda=(1,2)$. 
We have $x^{(0)}=(x_0, x_1)$, $x^{(1)}=(x_2)$, and the functional equations for this
case are:
\begin{align*}
G=& 1+zuGF_0 f_{0,1} \\
f_{k,0}=&g_k+g_{k+1}(x_0+x_1)+g_{k+2}x_0x_1 \\
f_{k,1}=&g_k+g_{k+1}x_2
\end{align*}
\vskip 1 cm
Now we make the substitution $x=uf_{0,1}$. Then the functional equations above become 
\begin{align*}
G=& 1+zxGF_0 \\
f_{k,0}=&g_k+g_{k+1}(x_0+x_1)+g_{k+2}x_0x_1
\end{align*}

By Theorem \ref{solution2}, the solution to these functional equations is:
\[G=\frac {c(t)E(y)}{E(y,c(t))},\]
where $y=(y_0,y_1)$, $E(y)=(1+y_0)(1+y_1)$ and $y$ is uniquely determined by $xx_i=\dfrac {y_i
(1+y_i)}{E^2(y)}$ for $i=0,1$, and $t=zxE^2(y)$.

Also, from Remark \ref{g1ex}, we have
\begin{align}
g_0=&1, \notag \\
g_1=& x E(y) \sum_{n=0}^{\infty} (1-n)e_n(y) \label{g1withx}\\
=&x(1+y_0)(1+y_1) (1-y_0y_1)\notag\\
=&uf_{0,1}(1+y_0)(1+y_1) (1-y_0y_1).\notag
\end{align}

Therefore, 
\begin{align*}
f_{0,1}=&g_0+g_1x_2\\
=&1+uf_{0,1}(1+y_0)(1+y_1) (1-y_0y_1) x_2.
\end{align*}
By solving for $f_{0,1}$, we get
$$f_{0,1}=\frac {1}{1-u(1+y_0)(1+y_1) (1-y_0y_1) x_2}.$$
So, by (\ref{g1withx}),
$$g_1=\frac {u(1+y_0)(1+y_1) (1-y_0y_1)}{1-u(1+y_0)(1+y_1) (1-y_0y_1) x_2}.$$

Now, since $xx_i=\dfrac {y_i (1+y_i)}{E^2(y)}$ for $i=0,1$ and $x=uf_{0,1}$,
we have 
$$ux_i=\dfrac {y_i (1+y_i)(1-u(1+y_0)(1+y_1) (1-y_0y_1) x_2)}{E^2(y)} \text{ for } i=0,1.$$ 
Since $(y_0,y_1)$ is uniquely determined by the equation above,
we can solve for $g_1$ in terms of $u, x_0, x_1$ and $x_2$ by using Lagrange inversion.

\vskip 2cm

We can also approach the problem without counting descents in the same way as what we did in
Section
\ref{another approach}.
If we let $P$ be the generating function of two-stack-sortable $k$-tuple $r$-permutations
under $S_{\lambda}$ for any $k$, with $u$ and ${\bar z}$ keeping track of the number of different letters and
the number of components, then
\begin{Lem}
P satisfies
\begin{equation} \label{FE51}
P= \frac {1}{1-\bar z}+\frac{uP}{{\bar z}^{\lambda_0}}
\left(P-\sum_{i=0}^{\lambda_0}p_i {\bar z}^i \right) \prod_{i=1}^{l} p_{\lambda_i-\lambda_{i-1}}.
\end{equation}
\end{Lem}

To solve this functional equation, first we make the substitution 
$$x=u\prod_{i=1}^{l} p_{\lambda_i-\lambda_{i-1}}.$$ 
Then (\ref{FE51}) becomes 
$$P= \frac {1}{1-\bar z}+\frac{xP}{{\bar z}^{\lambda_0}}
\left(P-\sum_{i=0}^{\lambda_0}p_i {\bar z}^i \right). $$

From Theorem \ref{solution2}, we know that
\begin{Thm}
The solution to the functional equation above is
$$
P=\left( 1+ \frac {t}{y(1+y)}\right) \frac {c(t)(1+y)^{\lambda_0+1}} { (1+
yc(t))^{\lambda_0+1}},
$$
where $y$ is uniquely determined as a formal power series in $x$ by 
$$ 
x=\frac {y}{(1+y)^{2 \lambda_0+1}},
$$ 
\begin{align*}
t=&x(1+y)^{2(\lambda_0+1)}\dfrac{\bar z}{1-\bar z} \\
=& y(1+y)\dfrac{\bar z}{1-\bar z} 
\end{align*}
and
$$x=u\prod_{i=1}^{l} p_{\lambda_i-\lambda_{i-1}}.
$$
\end{Thm}

Since $p_k$ is a polynomial in $y$ (Theorem \ref{pkpoly}) and 
$$u=\frac {y}{(1+y)^{2\lambda_0+1}\prod_{i=1}^{l}
p_{\lambda_i-\lambda_{i-1}}},
$$
we can solve for $p_k$ in terms of $u$.

Again, the simplest nontrivial case is when $l=1$ and $\lambda=(1,2)$, and the functional
equation for this case is:
$$
P= \frac {1}{1-\bar z} +up_1 \frac {P(P-p_0-p_1 \bar z)}{\bar z}.
$$
Now make the substitution $x=up_1$. Then the functional equation above becomes 
$$
P= \frac {1}{1-\bar z} +x \frac {P(P-p_0-p_1 \bar z)}{\bar z}.
$$
From Corollary \ref{p_1}, $p_1=1+y-y^2$, where $y$ is uniquely determined by 
$x=\dfrac{y}{(1+y)^3}$, or
$$ u = \frac {y}{(1+y)^3(1+y-y^2)}.$$ 

Now we can use Lagrange inversion [\ref{Goulden}, p.~ 17; \ref{Ira}] to solve for $p_1$. 

Since $y=u(1+y)^3(1+y-y^2)$,
\begin{align*}
[u^n]p_1=&\frac{1}{n} [y^{n-1}] \frac {d p_1}{d y} (1+y)^{3n}(1+y-y^2)^n \\
=&\frac{1}{n} [y^{n-1}] (1-2y) (1+y)^{3n} \sum_{i=0}^n (-1)^i \binom{n}{i}(1+y)^i y^{2(n-i)} \\
=& \frac{1}{n} [y^{n-1}] (1-2y) \sum_{i=0}^n \sum_{j=0}^{3n+i} 
(-1)^i\binom{n}{i} \binom{3n+i}{j} y^j y^{2(n-i)} \\
=&\sum_{i=0}^n (-1)^{n-i} \frac{1}{n} \binom{n}{i}
\left[ \binom{3n+i}{2i-1-n} -2 \binom{3n+i}{2i-2-n} \right] .
\end{align*}

\section{A characterization of $t$-stack-sortable permutations}

 It is natural to
consider counting $k$-stack-sortable permutations for $k>2$ now. But although people have been
trying, little has been found yet, even for three-stack-sortable permutations.

Similar to West's characterization of two-stack-sortable permutations (which did
not lead to an enumeration) [\ref{West2}, \ref{West}], we can give a
characterization for $t$-stack-sortable permutations.

For a sequence $a_1 a_2 \cdots a_n$ of different letters, define $(a_i,a_j)$ to be an inversion
of the sequence if $1<i<j<n$ and $a_i > a_j$. We also define
$\rank(a_i)$ as follows:
$\rank(a_i)=m$ if there are exactly $m-1$ letters among $a_1, a_2, \dots , a_n$ that are smaller
than $a_i$.

\begin{Thm} \label{tsschar}
A permutation $\pi$ is $t$-stack-sortable if and only if it does not contain a subsequence
$\pi^{\prime}=a_1 a_2
\cdots a_{t+2}$ which satisfies the following conditions: \\
(1) $\rank(a_{t+2})=1$; \\
(2) $\rank(a_{t+1})=t+2$; \\
(3) For any $i$ and $j$ such that $1 \le i < j \le
t$ and $(a_i, a_j)$ is an inversion, there does not exist a subsequence $\pi^{\prime
\prime}=c_1 c_2
\cdots c_s$ of $\pi_{\prime}$ where
$a_i< c_1 < c_2<\cdots <c_s$ and $s=t+2-\rank(a_i)$, such that $\pi^{\prime \prime}$ appears
between
$a_i$ and $a_j$ in $\pi$. \\
\end{Thm}

To prove this theorem, we need some lemmas from West [\ref{West2}, \ref{West}].

\begin{Lem} (West) \label{char1}
If $\pi$ is a permutation of $[n]$ and $1 \le a < b \le n$, and if $a$ precedes $b$ in $\pi$, then
$a$ precedes $b$ in $S(\pi)$.
\end{Lem}

\begin{Lem} (West) \label{char2}
If $\pi$ is a permutation of $[n]$ and $1 \le a < b \le n$, and if $b$ precedes $a$ in $\pi$, then
$b$ precedes $a$ in $S(\pi)$ if there exists $c>b$ such that $b$ precedes $c$ and $c$ precedes $a$
in $\pi$. If there is no such $c$, then $a$ precedes $b$ in $S(\pi)$.
\end{Lem}

\begin{Lem} (West) \label{char3}
If $b$ and $a$ form an inversion in $S(\pi)$, then there exist $c>b$ such that $b$ precedes $c$
and $c$ precedes $a$ in $\pi$.
\end{Lem}
Now let's prove Theorem \ref{tsschar}.
\begin{proof}
Tarjan and West proved that  the theorem is correct when $t=1$ and $t=2$ [\ref{tarjan}, \ref{West2}, \ref{West}].
We now prove the theorem by induction.

Suppose a permutation $\pi$ contains a subsequence $\pi^{\prime}=a_1 a_2
\cdots a_{t+2}$ satisfying the conditions (1) -- (3). Let $\pi_1=S(\pi)$ and suppose that
$b_1 b_2 \cdots b_{t+2}$ is the subsequence of $\pi_1$ that was $\pi^{\prime}$ in
$\pi$. 

Since $\pi^{\prime}$ satisfies conditions (1) -- (2), by Lemma \ref{char1}, $a_i$ precedes
$a_{t+1}$ in $S(\pi)$ for $1 \le i \le t$; by Lemma \ref{char2}, $a_i$ precedes $a_{t+2}$ in
$S(\pi)$ for $1 \le i \le t$. Thus, we get
either
$\rank(b_{t+2})=t+2$,
$\rank(b_{t+1})=1$ or 
$\rank(b_{t+2})=1$, $\rank(b_{t+1})=t+2$. Let $\pi_1^{\prime}=b_1 b_2 \cdots b_t b_{t+1}$ by
supposing 
$\rank(b_{t+1})=1$ (If $\rank(b_{t+2})=1$, we let $\pi_1^{\prime}=b_1 b_2 \cdots b_t b_{t+2}$).
Therefore,
$\pi_1^{\prime}$ satisfies condition (1) with $t$ replaced by $t-1$.

Let $\rank(a_{i_0})=t+1$. 
Then $a_{i_0} > a_i$ for $1 \le i \le t$ and $i \ne i_0$. By Lemma \ref{char1}, $a_i$ precedes
$a_{i_0}$  in $\pi_1$ for $1 \le i < i_0$. If $i_0
\ne t$, then
$a_{i_0}$ and $a_i$ form an inversion in $\pi^{\prime}$ for $i_0 < i \le t$. Since
$t+2-\rank(a_{i_0})=1$, by condition (3), there is no letter bigger than $a_{i_0}$ appears between
$a_{i_0}$ and
$a_i$ in $\pi$. By Lemma \ref{char2},
$a_i$ precedes $a_{i_0}$ in $\pi_1$ for $i_0 <i \le t$. Therefore, $a_{i}$ precedes $a_{i_0}$
in $\pi_1$ for
$1 \le i \le t$ and $i \ne i_0$.
So $b_t=a_{i_0}$ and $\rank(b_t)=t+1$. Hence $\pi_1^{\prime}$ satisfies condition (2).

Suppose that $\pi_1^{\prime}$ does not satisfy condition (3); that is, 
there exist $1 \le i < j \le t-1$ such that $b_i > b_j$, and there exists a subsequence
$\pi_1^{\prime \prime}=c_1 c_2
\cdots c_s$ where
$b_i< c_1 < c_2<\cdots <c_s$ and $s=t+1-\rank(b_i)$, such that $\pi_1^{\prime \prime}$ appears
between
$b_i$ and $b_j$ in $\pi_1$. Then $c_s$ and $b_j$ form an inversion in $\pi_1$. By
Lemma
\ref{char3}, there exists some $c > c_s$ such that $c_s$ precedes $c$ and $c$ precedes $b_j$ in
$\pi$. This contradicts the fact that $\pi^{\prime}$ satisfies condition (3).

Now, since $\pi_1^{\prime}$ satisfies condition (1) -- (3), by the induction hypothesis, $\pi_1$ is
not
$(t-1)$-stack-sortable. So $\pi$ is not $t$-stack-sortable.

Conversely, we can show that if $\pi$ fails to be $t$-stack-sortable, then it contains a
subsequence that satisfies the three conditions.

If $\pi$ is not $t$-stack-sortable, then $\pi_1=S(\pi)$ is not $(t-1)$-stack-sortable. By the
induction hypothesis, $\pi_1$ contains a subsequence $\pi_1^{\prime}=b_1 b_2 \cdots b_{t+1}$ which
satisfies the following conditions: \\
(1) $\rank(b_{t+1})=1$; \\
(2) $\rank(b_t)=t+1$;\\
(3) For any $i$ and $j$ such that $1 \le i < j \le
t-1$ and $b_i > b_j$, there does not exist a subsequence $\pi_1^{\prime \prime}=c_1 c_2
\dots c_s$ where
$b_i< c_1 < c_2<\cdots <c_s$ and $s=t+1-\rank(b_i)$, such that $\pi_1^{\prime \prime}$ appears
between
$b_i$ and $b_j$ in $\pi_1$.

Notice that $b_i$ and $b_{t+1}$ form an inversion in $\pi_1$ for $1 \le i \le t$. By Lemma
\ref{char3}, $b_i$ precedes $b_{t+1}$ in $\pi$ for $1 \le i \le t$. In particular, there exists
some $a > b_t$ such that $b_i$ precedes $a$ and $a$ precedes $b_{t+1}$ in $\pi$ for 
$1 \le i \le t$. Suppose  that
$a_1 a_2 \dots a_{t+1}$ is the subsequence of $\pi$ that gets tranformed to $\pi_1^{\prime}$. 
Then it is clear that $a_{t+1}=b_{t+1}$. Let $\pi^{\prime}=a_1a_2 \dots a_{t} a a_{t+1}$. Then  
$\pi^{\prime}$ satisfies conditions (1) and (2).

If $\pi^{\prime}$ does not satisfy condition (3); that is, 
there exist $1 \le i < j \le t$ such that $a_i > a_j$, and there exists a subsequence
$\pi^{\prime \prime}=c_1 c_2
\dots c_s$ where
$a_i< c_1 < c_2<\cdots <c_s$ and $s=t+2-\rank(a_i)$, such that $\pi^{\prime \prime}$ appears
between
$a_i$ and $a_j$, then by Lemma \ref{char1}, the subsequence $c_1 c_2 \dots c_{s-1}$ appears
between $a_i$ and $a_j$ in $\pi_1$. This contradicts the fact that $\pi_1^{\prime}$ satisfies
condition (3). Therefore $\pi^{\prime}$ satisfies condition (3).

The theorem is proved.
\end{proof}

\begin{Lem} (Tarjan)
A permutation $\pi$ is stack-sortable if and only if $\pi$ contains no subsequence of type 231.
\end{Lem}
\begin{proof}
This lemma is the case $t=1$ of Theorem \ref{tsschar}.

By Theorem \ref{tsschar}, a permutation $\pi$ is stack-sortable if and only if it does not contain
a subsequence
$\pi^{\prime}=a_1a_2a_3$, which satisfies the conditions that $\rank(a_3)=1$ and $\rank(a_2)=3$,
which means that
$a_3 <a_1 <a_2$. Therefore,  $\pi^{\prime}$ is a type 231 subsequence.  
\end{proof}

\begin{Lem} (West)
A permutation fails to be two-stack-sortable if it contains a subsequence of type 2341 or a
subsequence of type 3241 which is not part of a subsequence of type 35241. If it contains no such
subsequence, $\pi$ is two-stack-sortable.
\end{Lem}
\begin{proof}
This lemma is the case $t=2$ of Theorem \ref{tsschar}.

By Theorem \ref{tsschar}, a permutation $\pi$ is two-stack-sortable if and only if it does not
contain a subsequence $\pi^{\prime}=a_1a_2a_3a_4$, which satisfies the conditions that
(1) $\rank(a_4)=1$, (2) $\rank(a_3)=4$ and, (3) $a_1 < a_2$ or $ a_1 >a_2 $. If $a_1 < a_2$, then 
$\pi^{\prime}$ is a subsequence of type 2341. If $a_1 > a_2$, then $\pi^{\prime}$ is a subsequence of type
3241. By condition (3), there does not exist a $c>a_1$ appearing between $ a_1$ and $a_2$. 
If there does exist such a $c$, then when $c>a_3$, $a_1 c a_2 a_3 a_4$ is a subsequence of type 35241; when 
$c < a_3$, $c a_2 a_3 a_4$ is a subsequence of 3241, which needs to be considered again in the same way. 
Therefore, these
three conditions suggest that
$\pi^{\prime}$ is a either type 2341 or type 3241 subsequence. If it is a type 3241 subsequence,
then it is not part of a subsequence of type 35241.
\end{proof}

To enumerate the three-stack-sortable permutations by the same decomposition we
used for two-stack-sortable permutations, we consider the following object, called a
$\mu$-tuple permutation, where $\mu=(\mu_1,\dots,\mu_k)$ and 
$\mu_i$ are positive integers.

\begin{Def}
$((\alpha_1,\dots,\alpha_{\mu_1}),
(\alpha_{\mu_1+1},\dots,\alpha_{\mu_1+\mu_2}),
\dots,
(\alpha_{\mu_1+\cdots+\mu_{k-1}+1},\dots,
\alpha_{\mu_1+\cdots+\mu_{k}}))$ is a $\mu$-tuple permutation of $[n]$
if
$\alpha_1 \cdots \alpha_{\mu_1+\cdots+\mu_k}$ is a permutation of $[n]$. A
$\mu$-tuple permutation is three-stack-sortable if 
$$S(S(S(\alpha_1) \cdots S(\alpha_{\mu_1}))
S(S(\alpha_{\mu_1+1}) \cdots S(\alpha_{\mu_1+\mu_2}))
\cdots
S(S(\alpha_{\mu_1+\cdots+\mu_{k-1}+1}) \cdots
S(\alpha_{\mu_1+\cdots+\mu_{k}})))=I.$$
\end{Def}

By the same reasoning as in the case of two-stack-sortable permutations, if
$g_{\mu}$ is the generating function for three-stack-sortable $\mu$-tuple
permutations, we get 

$$g_{\mu}=1+x\left(
\sum_{i=0}^{\mu_k-1} g_{(\mu_1,\dots,\mu_k-i+1,i)} +
\sum_{i=0}^{\mu_{k-1}-1}
g_{(\mu_1,\dots,\mu_{k-1}-i+1,i)}g_{(\mu_k)} + \cdots
+\sum_{i=0}^{\mu_1-1} g_{(\mu_1-i+1)}g_{(\mu_2,\dots,\mu_k)}
\right). $$

It seems that this functional equation is hard to solve.

\section{A modification on the functional equations.} 
If we replace the
elementary symmetric functions 
$e_i(x)$ by complete homogeneous symmetric functions $h_i(x)$, then the functional
equations in Theorem \ref{functional equations} will become:
\begin{align*} 
f_k=&\sum_{j=0}^{\infty} g_{k+j} h_j(x) \\
G=&1+zFG
\end{align*}
Combinatorially, in terms of the tree representation, instead of having at most one child of each type,
now we can have any number of children of each type. 
To solve this,
everything else will be the same except one modification on the substitution
$$x_i=\frac {y_i(1+y_i)} {E^2(y)},$$
which in this case will be 
$$x_i=\frac {y_i(1-y_i)} {H^2(y)},$$
where $H(y)$ is defined analogously to $E(y)$. Then analogous to the identity
in Lemma \ref{identity}, we get
$$ H\left(x, \frac{1}{z}\right) = H(y,c(t))H\left( y, \frac{1}{tc(t)}\right)$$

The solutions to these functional equations
will be:
\begin{align*}
G=&\frac {c(t)H(y)}{H(y,c(t))}, \\
F=&\frac{H(y)}{t}\left(H(y)-\frac{H(y,c(t))}{c(t)}\right).
\end{align*}

Therefore, 
\begin{align*}
g_1=&f_0 \\
=&H(y) \sum_{n=0}^{\infty}(1-n)h_n(y)
\end{align*}

Again, we use multivariable Lagrange inversion to get the coefficient of
$x^k$ in
$g_1$, where $x=(x_0,\dots,x_r)$ and $k=(k_0,\dots, k_r)$ and it is  
\[
\frac{1}{n^2} \prod_{k=0}^r \frac{n}{n+k_i} \binom {2n+2k_i}{k_i},
\]
where $n=1+k_0+k_1+\cdots+k_r$.

\vs 5
\begin{Remark}
If we let 
\[ A(n; k_0, \dots, k_r) = \frac{1}{n^2} \prod_{k=0}^r \frac{n}{n-k_i} \binom {2n-1-k_i}{k_i} \]
and 
\[ B(n; k_0, \dots, k_r) = \frac{1}{n^2} \prod_{k=0}^r \frac{n}{n+k_i} \binom {2n+2k_i}{k_i} ,\]
then
\[ A(-n;  k_0, \dots, k_r) = (-1)^{n-1} B(n ;  k_0, \dots, k_r). \]
\end{Remark}

\vs 9

\vs 9

\end{document}